\documentclass[12pt]{article}
\usepackage[latin9]{inputenc}
\usepackage{float}
\usepackage{bbding}
\usepackage{amsmath}
\usepackage{amssymb}
\usepackage{graphicx}

\makeatletter

\usepackage{amsfonts}
\usepackage{graphicx}
\providecommand{\U}[1]{\protect\rule{.1in}{.1in}}
\newtheorem{theorem}{Theorem}\newtheorem{algorithm}[theorem]{Algorithm}\newtheorem{case}[theorem]{Case}\newtheorem{definition}[theorem]{Definition}\newtheorem{problem}[theorem]{Problem}

\makeatother

\begin{document}
\title{\textbf{Superiorization: The asymmetric roles of feasibility-seeking
and objective function reduction}}
\author{Yair Censor\medskip{}
\\
{Department of Mathematics, University of Haifa}\\
Mt. Carmel, Haifa 3498838, Israel\\
yair@math.haifa.ac.il}
\date{October 24, 2022. Revised: December 29, 2022.}

\maketitle
 
\begin{abstract}
The superiorization methodology can be thought of as lying conceptually
between feasibility-seeking and constrained minimization. It is not
trying to solve the full-fledged constrained minimization problem
composed from the modeling constraints and the chosen objective function.
Rather, the task is to find a feasible point which is ``superior''
(in a well-defined manner) with respect to the objective function,
to one returned by a feasibility-seeking only algorithm. We telegraphically
review the superiorization methodology and where it stands today and
propose a rigorous formulation of its, yet only partially resolved,
guarantee problem. The real-world situation in an application field
is commonly represented by constraints defined by the modeling process
and the data, obtained from measurements or otherwise dictated by
the model-user. The feasibility-seeking problem requires to find a
point in the intersection of all constraints without using any objective
function to aim at any specific feasible point.

At the heart of the superiorization methodology lies the modeler desire
to use an objective function, that is exogenous to the constraints,
in order to seek a feasible solution that will have lower (not necessarily
minimal) objective function value. This aim is less demanding than
full-fledged constrained minimization but more demanding than plain
feasibility-seeking. 

Putting emphasis on the need to satisfy the constraints, because they
represent the real-world situation, one recognizes the ``asymmetric
roles of feasibility-seeking and objective function reduction'',
namely, that fulfilling the constraints is the main task while reduction
of the exogenous objective function plays only a secondary role. There
are two research directions in the superiorization methodology that
nourish from this same general principle: Weak superiorization and
strong superiorization. 

Since its inception in 2007, the superiorization methodology has evolved
and gained ground, as can be seen from the, compiled and continuously
updated, bibliography at: http://math.haifa.ac.il/yair/bib-superiorization-censor.html\#top.\bigskip{}
\end{abstract}

\tableofcontents{}

\section{Introduction\label{sec:Introduction}}

This paper is a tutorial/review but it does not follow the classical
sense of these notions. It mostly rather collects pieces from the
published literature creating a collage\footnote{As such, it inevitably includes some portions from previous publications
which are all clearly acknowledged.}. This collage brings together definitions and results about the \textbf{superiorization
methodology}, directs the reader to the existing literature, says
a word about the history, and proposes a rigorous formulation of its,
yet only partially resolved, guarantee problem.

It is a tribute to Simeon Reich, friend, colleague and collaborator, 
whose papers (with Butnariu and Zaslavski) \cite{brz06,brz08} paved
the way for the superiorization methodology, see Section \ref{sec:sup-framework}
below, on the occasion of his 75th birthday.

The fundamental underlying question considered here is what to do
with \textbf{input} that consists of a \textbf{constraints set} and
an \textbf{objective function}?\textbf{ }Do constrained minimization
or do superiorization? 

There are many available routes toward doing constrained minimization.
These include, speaking generally, algorithms based on the use of
Lagrange multipliers or on the regularization approach. In the latter
approach a \textbf{proximity function} that measures constraints violation
is appended to the objective function, along with a regularization
parameter, and algorithms for unconstrained minimization are applied
to the resulting ``regularized objective function''. 

Other methods, like penalty methods and barrier methods, require a
feasible initialization point (i.e., a point in the constraints set)
and perform searches for the constrained minimum within the constraints
set by ``preventing'' the iterates of these searches from ``leaving''
the constraints set.

Recognizing that the two pieces of the input (the constraints set
and the objective function) are independent, there is the approach
of handling the constrained minimization problem as an unconstrained
bi-objective (i.e., multi-objective with two objectives) problem for
the original objective function and for a proximity function.

In the regularization approach, mentioned above, the user decides
how much importance should be attached to fulfilling the constraints
versus minimizing the original objective function by defining the,
so called, \textbf{regularization parameter}, which determines how
much weight is given to each of the two tasks.

The superiorization methodology (SM) offers a different approach than
those mentioned above. Putting emphasis on the need to satisfy the
constraints, since they represent the real-world situation, the SM
recognizes the asymmetric roles of feasibility-seeking and objective
function reduction. Namely, that the feasibility-seeking to fulfill
the constraints, is the main task while reduction of the exogenous
objective function plays only a secondary role. 

In a nutshell, the SM does this by taking a convergent feasibility-seeking
algorithm and perturbing its iterates so that the objective function
values get reduced while retaining the overall convergence to a feasible
point. The question whether such ``local'' objective function reduction
steps, that are interlaced into the feasibility-seeking algorithm,
accumulate toward a global objective function reduction, not necessarily
minimization, is the gist of the yet not fully solved guarantee problem
of the SM. 

In the SM one takes an iterative algorithm, investigates its \textbf{perturbation
resilience}, and then, uses proactively such permitted perturbations,
to force the perturbed algorithm to do something useful in addition
to what it is originally designed to do. The original unperturbed
algorithm is called the \textbf{basic algorithm} and the perturbed
algorithm is called the \textbf{superiorized version of the basic
algorithm}.

If the basic algorithm is computationally efficient and useful in
terms of the application at hand, and if the perturbations are simple
and not expensive to calculate, then the advantage of this methodology
is that, for essentially the computational cost of the basic algorithm,
we are able to get something more by steering its iterates according
to the perturbations.

This is a very general principle, which has been successfully used
in a variety of important practical applications and awaits to be
implemented and tested in additional fields; see the recent papers
in the, compiled and continuously updated, bibliography of scientific
publications on the superiorization methodology and perturbation resilience
of algorithms \cite{SM-bib-page}.

Although not limited to this case, an important special case of the
superiorization methodology is when the basic algorithm is a feasibility-seeking
algorithm for a family of constraints and the perturbations that are
interlaced into it aim at reducing, not necessarily minimizing, a
given objective function. This case is what drives the intuition behind
the superiorization methodology and is the subject matter of this
paper\footnote{Support for this reasoning may be borrowed from the American scientist
and Noble-laureate Herbert Simon who was in favor of ``satisficing''
rather than ``maximizing''. Satisficing is a decision-making strategy
that aims for a satisfactory or adequate result, rather than the optimal
solution. This is because aiming for the optimal solution may necessitate
needless expenditure of time, energy and resources. The term ``satisfice''
was coined by Herbert Simon in 1956 \cite{simon}, see also: https://en.wikipedia.org/wiki/Satisficing.}.

Superiorization has a world-view that is quite different from that
of classical constrained optimization. Both in superiorization and
in classical constrained optimization there is an assumed domain $\Omega\subseteq R^{n}$
in the $n$-dimensional Euclidean\footnote{We limit most of our discussion to the finite-dimensional Euclidean
space in order to present as clearly as possible our points although
the notions and ideas can be carried over to Hilbert space or other
infinite-dimensional spaces.} space and an exogenous objective (a.k.a. merit, criterion, target\footnote{We use ``objective function'' and ``target function'' interchangeably
in the sequel.}, etc.) function $\phi:R^{n}\rightarrow R$ that maps $\varOmega$
into $R$. In classical optimization it is assumed that there is a
constraints set $C\subset\varOmega$ and the task is to find an $x\in C$
for which $\phi(x)$ is minimal over $C$. In superiorization the
task is different: to find a feasible point in $C$ which is ``superior''
(in a well-defined manner) with respect to the objective function,
to one returned by a feasibility-seeking only algorithm.

There are two research directions in the superiorization methodology
that nourish from the same general principle. One is the direction
when the constraints are assumed to be consistent, $C\neq\emptyset$
and the notion of \textbf{bounded perturbation resilience} is used.
In this case one treats the superiorized version of the basic algorithm
as a recursion formula, without a stopping rule, that produces an
infinite sequence of iterates and asymptotic convergence questions
are in the focus of study.

The second direction does not assume consistency of the constraints
but uses instead a \textbf{proximity function} that measures the violation
of the constraints. Instead of seeking asymptotic feasibility, it
looks at \textbf{$\varepsilon$-compatibility} and uses the notion
of \textbf{strong perturbation resilience}. The same core superiorized
version of the basic algorithm might be investigated in each of these
directions, but the second is apparently more practical since it relates
better to problems formulated and treated in practice. We use the
terms \textbf{weak superiorization} and \textbf{strong superiorization}
as a nomenclature for the first and second directions, respectively\footnote{These terms were proposed in \cite{cz14-feje}, following a private
discussion with our colleague and coworker in this field Gabor Herman.}.

\section{Constraints oriented modeling\label{sect:con-orient-mod}}

Let $R^{n}$ be the $n$-dimensional Euclidean space where $x\in R^{n}$
is represented by its components $x=(x_{j})_{j=1}^{n}.$ A \textbf{constraint}
is a condition that restricts a vector $x$ to belong to given sets
that are represented by functional inequalities. \textbf{Modeling
}is the process of representing a real-world problem, in some field
of application, in a mathematical language amenable to mathematical
analysis and to the development of tractable algorithmic solutions.\textbf{ }

\textbf{Constraints oriented modeling }is a modeling process that
represents the real-world problem by a system of constraints. The
``individual'' constraints $C_{i}\subset R^{n},$ for all $i=1,2,\ldots,m,$
are
\begin{equation}
C_{i}:=\{x\in R^{n}\mid q_{i}(x)\leq\gamma_{i}\}
\end{equation}
 where the function $q_{i}:R^{n}\rightarrow R$ is a mapping into
the reals and $\gamma_{i}\in R$ and $C_{i}$ is the $\gamma_{i}$-level
set of $q_{i}.$ This gives rise to the \textbf{feasibility-seeking
problem}
\begin{equation}
\textup{find}\;x\in C:=\cap_{i=1}^{m}C_{i}.\label{eq:CFP}
\end{equation}
We assume that this feasibility-seeking problem is \textbf{consistent},
i.e., $C\neq\emptyset.$ In the inconsistent case the statement $x\in C$
is meaningless and alternative solution concepts must be used, but
we do not wander in this direction here, see, e.g., \cite{cen-zak-inconsist-2018}
and references therein. 

The word ``find'' in (\ref{eq:CFP}) is meant as either actually
presenting such an $x\in C$ or generating an infinite sequence $\{x^{k}\}{}_{k=0}^{\infty}$,
with $x^{k}\in R^{n}$ for all $k\geq0,$ that asymptotically converges
such that $\lim_{k\rightarrow\infty}x^{k}=x^{*}\in C.$ In strong
superiorization, however, ``find'' has another meaning, see Section
\ref{sect:strongSM} below.

A particular instance is when all constraints are convex sets represented
as level sets of convex functions. In that case (\ref{eq:CFP}) is
the well-known \textbf{convex feasibility problem} (CFP), consult
\cite{bb96} and \cite[Subsection 1.3.4]{cegielski-book-2012}, which
is, in turn, a special instance of the \textbf{common fixed point
problem} (CFPP) of a family of operators, consult \cite[Subsection 4.6]{BC-book-2nd ed-2017}.

In the constraints oriented modeling process the constraints functions
$\{q_{i}\}{}_{i=1}^{m}$ are usually decided upon by the modeler to
best represent the nature of the problem that is modeled. The \textbf{level
sets parameters} $\{\gamma_{i}\}{}_{i=1}^{m}$ are commonly obtained
through physical measurements or via prescriptions defined by the
user. Thus, $\{q_{i}\}{}_{i=1}^{m}$ and $\{\gamma_{i}\}{}_{i=1}^{m}$
are the \textbf{input} for the feasibility-seeking problem.

A question that might arise at this point is why not translate the
feasibility-seeking problem upfront to an unconstrained minimization
of a proximity function that measures the violation of the constraints
or by using the indicator functions of the individual sets? Or why
not translate the feasibility-seeking problem upfront to a constrained
minimization of an objective function which is a constant over the
constraints at hand? While these approaches are legitimate they do
not necessarily lead to the same ``algorithmic territory''. It is
not clear, and to our understanding doubtful, whether better algorithms
for the feasibility-seeking problem can be discovered when using such
translational approaches. 

On the contrary, quite a few well-known algorithms for feasibility-seeking
were, and still are, discovered, developed and studied by using tools
from outside mathematical optimization. for example, the feasibility-seeking
problem can be formulated as a special case of the common fixed point
problem, mentioned above, for the case when the operators are projections,
and as such benefit from the large body of knowledge in the field
of fixed point theory. The excellent book of Cegielski \cite{cegielski-book-2012}
attests to this line of work.

\section{Feasibility-seeking algorithms}

We look at the \textbf{convex feasibility problem} (CFP) which is
to find a feasible point $x^{\ast}\in C$ when all sets $C_{i}$ are
convex and commonly also closed. This prototypical problem underlies
the modeling of a variety of real-world problems in many fields, see,
e.g., the pointers and references in Bauschke and Borwein \cite[Section 1]{bb96}
and in Cegielski's book \cite[Section 1.3]{cegielski-book-2012}.

If $C\neq\emptyset$ does not hold then the CFP is inconsistent and
a feasible point does not exist, see, e.g., the review of inconsistent
feasibility problems \cite{cen-zak-inconsist-2018}. However, algorithmic
research of inconsistent CFPs exists and is mainly focused on two
directions. One is oriented toward defining solution concepts other
than $x^{\ast}\in C$ that will apply, such as proximity function
minimization wherein a \textbf{proximity function} measures in some
way the total violation of all constraints, see, e.g., \cite[pp. 28-29.]{cegielski-book-2012}.
The second direction investigates the behavior of algorithms that
are designed to solve a consistent CFP when applied to inconsistent
problems. The latter is fueled by situations wherein one lacks a priori
information about the consistency or inconsistency of the CFP or does
not wish to invest computational resources to get hold of such knowledge
prior to running his algorithm.

\textbf{Projection methods}. Projections onto sets are used in a wide
variety of methods in optimization theory but not every method that
uses projections really belongs to the class of projection methods
as we mean it here. Here \textbf{projection methods} are iterative
algorithms\footnote{As common, we use the terms algorithm or algorithmic structure for
the iterative processes studied here although no termination criteria,
which are by definition necessary in an algorithm, are present and
only the asymptotic behavior of these processes is studied. This does
not create any ambiguity because whether we consider an infinite iterative
process or an algorithm with a termination rule is always clear from
the context.} that use projections onto sets while relying on the general principle
that when a family of (usually closed and convex) sets is present
then projections (or approximate projections) onto the given individual
sets are easier to perform than projections onto other sets (intersections,
image sets under some transformation, etc.) that are derived from
the family of individual sets.

A projection algorithm reaches its goal, related to the whole family
of sets, by performing projections onto the individual sets. Projection
algorithms employ projections (or approximate projections) onto convex
sets in various ways. They may use different kinds of projections,
e.g., orthogonal (least Euclidean distance) projections, Bregman projections,
entropy projections, subgradient projections, intrepid projections,
valiant projections, Douglas-Rachford operators etc. and, sometimes,
even use different projections within the same algorithm. They serve
to solve a variety of problems which are either of the feasibility-seeking
or the constrained optimization types. They have different algorithmic
structures, of which some are particularly suitable for parallel computing,
and they demonstrate nice convergence properties and/or good initial
behavior patterns in some significant fields of applications.

Apart from theoretical interest, the main advantage of projection
methods, which makes them successful in real-world applications, is
computational. They commonly have the ability to handle huge-size
problems of dimensions beyond which other, more sophisticated currently
available, methods cease to be efficient. This is so because the building
bricks of a projection algorithm are the projections onto the individual
sets (assumed and actually easy to perform) and the algorithmic structures
are either sequential or simultaneous or in-between, such as in the
\textbf{block-iterative projection} (BIP) methods, see, e.g., \cite{combettes-BIP},
\cite{shoham-BIP}, \cite{ndh12}, or in the more recent \textbf{string-averaging
projection} (SAP) methods, see details and references in Subsection
\ref{sect:weakSM}. An advantage of projection methods is that they
work with initial data and do not require transformation of, or other
operations on, the sets describing the problem.

\section{The superiorization methodology\label{sec:sup-framework}}

Since its inception in 2007, the superiorization method has evolved
and gained ground. Quoting and distilling from earlier publications,
we review here the two directions of the superiorization methodology.
Recent review papers on the subject which could be read together with
this paper are Herman's \cite{gth-sup4IA} and \cite{herman-jano-2020}.
Unless otherwise stated, we restrict ourselves, for simplicity, to
the $n$-dimensional Euclidean space $R^{n}$ although some materials
below remain valid in Hilbert space.

Recent publications on the superiorization methodology (SM) are devoted
to either weak or strong superiorization, without yet using these
terms. They are \cite{bk13,bdhk07,cdh10,cz12,dhc09,rand-conmath,gh13,hd08,hgdc12,wenma13,ndh12,pscr10}.
many of the papers contain a detailed description of the SM, its motivation,
and an up-to-date review of SM-related previous work. 

The superiorization method was born when the terms and notions \textquotedblleft superiorization\textquotedblright{}
and \textquotedblleft perturbation resilience\textquotedblright ,
in the present context, first appeared in the 2009 paper \cite{dhc09}
which followed its 2007 forerunner by Butnariu et al. \cite{bdhk07}.
The ideas have some of their roots in the 2006 and 2008 papers of
Butnariu et al. \cite{brz06,brz08} where it was shown that if iterates
of a nonexpansive operator converge for any initial point, then its
inexact iterates with summable errors also converge. 

Bounded perturbation resilience of a parallel projection method was
observed as early as 2001 in \cite[Theorem 2]{combettes2001} (without
using this term). All these culminated in Ran Davidi\textquoteright s
2010 PhD dissertation \cite{davidi-thesis} and the many papers that
appeared since then and are cited in \cite{SM-bib-page}. The latter
is a Webpage dedicated to superiorization and perturbation resilience
of algorithms that contains a continuously updated bibliography on
the subject. This Webpage\footnote{http://math.haifa.ac.il/yair/bib-superiorization-censor.html\#top,
last updated on December 4, 2022 with 161 items.} is source for the wealth of work done in this field to date, including
two special issues of journals \cite{SI-inverse-prob-2017} and \cite{SI-JANO-2020}
dedicated to research of the SM. Recent work includes \cite{nikazad-super-2022},
\cite{hoseini-super-2020}, \cite{he-xu-2017,bonacker-2017,wang-2017}.
Interestingly, \cite{andersen-2014} notices some structural similarities
of the SM with incremental proximal gradient methods.

Let $T$ denote a mathematically-formulated problem, of any kind or
sort, with solution set $\Psi_{T}.$ The following cases immediately
come to mind although any $T$ and its $\Psi_{T}$ can potentially
be used.

\begin{case} \label{case:cfp}$T$ is a \textit{convex feasibility
problem} (CFP) of the form: find a vector $x^{\ast}\in\cap_{i=1}^{m}C_{i},$
where $C_{i}\subseteq R^{n}$ are closed convex subsets. In this case
$\Psi_{T}=\cap_{i=1}^{m}C_{i}.$ \end{case}

\begin{case} \label{case:2}$T$ is a constrained minimization problem:
$\mathrm{minimize}\left\{ f(x)\mid x\in\Phi\right\} $ of an objective
function $f$ over a feasible region $\Phi.$ In this case $\Psi_{T}=\{x^{\ast}\in\Phi\mid f(x^{\ast})\leq f(x)$
for all $x\in\Phi\}.$ \end{case}

The superiorization methodology is intended for function reduction
problems of the following form.

\begin{problem} \label{prob:sm}\textbf{The Function Reduction Problem}.
Let $\Psi_{T}\subseteq{R^{n}}$ be the solution set of some given
mathematically-formulated problem $T$ and let $\phi:R^{n}\rightarrow R$
be an objective function. Let $\mathcal{A}:R^{n}\rightarrow R^{n}$
be an algorithmic operator that defines an iterative \textbf{basic
algorithm} for the solution of $T$. Find a vector $x^{\ast}\in\Psi_{T}$
whose function $\phi$ value is lesser than that of a point in $\Psi_{T}$
that would have been reached by applying the Basic Algorithm for the
solution of problem $T.$ \end{problem}

As explained below, the superiorization methodology approaches this
problem by automatically generating from the basic algorithm its \textbf{superiorized
version of the basic algorithm}. The vector $x^{\ast}$ obtained from
the superiorized version of the basic algorithm need not be a minimizer
of $\phi$ over $\Psi_{T}.$ Another point to observe is that the
very problem formulation depends not only on the data $T,$ $\Psi_{T}$
and $\phi$ but also on the pair of algorithms -- the original unperturbed
basic algorithm, represented by $\mathcal{A},$ for the solution of
problem $T,$ and its superiorized version.

A fundamental difference between weak and strong superiorization lies
in the meaning attached to term ``solution of problem $T$'' in
Problem \ref{prob:sm}. In weak superiorization solving the problem
$T$ is understood as generating an infinite sequence $\{x^{k}\}_{k=0}^{\infty}$
that converges to a point $x^{\ast}\in\Psi_{T},$ thus, $\Psi_{T}$
must be nonempty. In strong superiorization solving the problem $T$
is understood as finding a point $x^{\ast}$ that is \textbf{$\varepsilon$-compatible
with $\Psi_{T},$} for some positive $\varepsilon,$ thus, nonemptiness
of $\Psi_{T}$ need not be assumed.

We concentrate in the next sections mainly on Case \ref{case:cfp}.
Superiorization work on Case \ref{case:2}, e.g., where $T$ is a
maximum likelihood optimization problem and $\Psi_{T}$ -- its solution
set, appears in \cite{gh13,wenma13,tie}.

\subsection{Weak superiorization\label{sect:weakSM}}

In weak superiorization the set $\Psi_{T}$ is assumed to be nonempty
and one treats the ``Superiorized Version of the Basic Algorithm''
as a recursion formula that produces an infinite sequence of iterates.
The SM strives to asymptotically find a point in $\Psi_{T}$ which
is superior,\textbf{ }i.e., has a lower, but not necessarily minimal,
value of the $\phi$ function, to one returned by the Basic Algorithm
that solves the original problem $T$ only.

This is done by first investigating the bounded perturbation resilience
of an available Basic Algorithm designed to solve efficiently the
original problem $T$ and then proactively using such permitted perturbations
to steer its iterates toward lower values of the $\phi$ objective
function while not loosing the overall convergence to a point in $\Psi_{T}$.

\begin{definition} \label{def:resilient}\textbf{Bounded perturbation
resilience (BPR)}. Let $\Gamma\subseteq R^{n}$ be a given nonempty
set. An algorithmic operator $\mathcal{A}:R^{n}\rightarrow R^{n}$
is said to be \textbf{bounded perturbations resilient with respect
to $\Gamma$}\emph{ }if the following is true: If a sequence $\{x^{k}\}_{k=0}^{\infty},$
generated by the iterative process $x^{k+1}=\mathcal{A}(x^{k}),$
for all $k\geq0,$ converges to a point in $\Gamma$ for all $x^{0}\in R^{n}$,
then any sequence $\{y^{k}\}_{k=0}^{\infty}$ of points in $R^{n}$
that is generated by $y^{k+1}=\mathcal{A}(y^{k}+\beta_{k}v^{k}),$
for all $k\geq0,$ also converges to a point in $\Gamma$ for all
$y^{0}\in R^{n}$ provided that, for all $k\geq0$, $\beta_{k}v^{k}$
are \textbf{bounded perturbations}, meaning that $\beta_{k}\geq0$
for all $k\geq0$ such that ${\displaystyle \sum\limits _{k=0}^{\infty}}\beta_{k}<\infty,$
and that the sequence $\{v^{k}\}_{k=0}^{\infty}$ is bounded. \end{definition}

Let $\phi:R^{n}\rightarrow R$ be a real-valued convex continuous
function and let $\partial\phi(z)$ be the subgradient set of $\phi$
at $z$ and, for simplicity of presentation, assume here that $\Gamma=R^{n}.$
In other specific cases care must be taken regarding how $\Gamma$
and $\Psi_{T}$ are related. The following superiorized version of
the basic algorithm\textbf{ }$\mathcal{A}$ is based on \cite[Algorithm 4.1]{cz14-feje}.

\begin{algorithm}
\label{alg:super-process}$\left.{}\right.$\textbf{Algorithm 1. Superiorized
version of the basic algorithm }$\mathcal{A}$\textbf{.}

\textbf{(0) Initialization}: Let $N$ be a natural number and let
$y^{0}\in R^{n}$ be an arbitrary user-chosen vector.

\textbf{(1)} \textbf{Iterative step}: Given a current iteration vector
$y^{k}$ pick an $N_{k}\in\{1,2,\dots,N\}$ and start an inner loop
of calculations as follows:

\textbf{(1.1) Inner loop initialization}: Define $y^{k,0}=y^{k}.$

\textbf{(1.2) Inner loop step: }Given $y^{k,n},$ as long as $n<N_{k},$
do as follows:

\textbf{(1.2.1) }Pick a $0<\beta_{k,n}\leq1$ in a way that guarantees
that
\begin{equation}
\sum_{k=0}^{\infty}\sum_{n=0}^{N_{k}-1}\beta_{k,n}<\infty.\label{eq:2.6}
\end{equation}

\textbf{(1.2.2)} Pick an ${\displaystyle s^{k,n}\in\partial\phi(y^{k,n})}$
and define $v^{k,n}$ as follows:
\begin{equation}
v^{k,n}=\left\{ \begin{array}{cc}
-\frac{{\displaystyle s^{k,n}}}{{\displaystyle \left\Vert s^{k,n}\right\Vert }}, & \text{if }0\notin\partial\phi(y^{k,n}),\\
0, & \text{if }0\in\partial\phi(y^{k,n}).
\end{array}\right.
\end{equation}

\textbf{(1.2.3) }Calculate the perturbed iterate
\begin{equation}
y^{k,n+1}=y^{k,n}+\beta_{k,n}v^{k,n}\label{eq:2.11}
\end{equation}
and if $n+1<N_{k}$ set $n\leftarrow n+1$ and go to \textbf{(1.2)},
otherwise go to \textbf{(1.3)}.

\textbf{(1.3) }Exit the inner loop with the vector $y^{k,N_{k}}$

\textbf{(1.4) }Calculate
\begin{equation}
y^{k+1}=\mathcal{A}(y^{k,N_{k}})\label{eq:2.12}
\end{equation}
set $k\leftarrow k+1$ and go back to \textbf{(1)}. 
\end{algorithm}

Considering Case \ref{case:cfp} in Section \ref{sec:sup-framework}
where $T$ is a convex feasibility problem, the Dynamic String-Averaging
Projection (DSAP) method of \cite{cz12} that we describe here constitutes
a family of algorithmic operators that can play the role of the above
$\mathcal{A}$ in a basic algorithm for the solution of the CFP $T$.

Let $\{C_{i}\}{}_{i=1}^{m}$ be nonempty closed convex subsets of
a Hilbert space $X$, set $C=\cap_{i=1}^{m}C_{i},$ and assume $C\neq\emptyset$.
For $i=1,2,\dots,m,$ denote by $P_{i}:=P_{C_{i}}$ the orthogonal
(least Euclidean distance) projection onto the set $C_{i}.$ An \textbf{index
vector} is a vector $t=(t_{1},t_{2},\dots,t_{q})$ such that $t_{s}\in\{1,2,\dots,m\}$
for all $s=1,2,\dots,q$, whose length is $\ell(t)=q.$ The product
of the individual projections onto the sets whose indices appear in
the index vector $t$ is $P[t]:=P_{t_{q}}\cdots P_{t_{2}}P_{t_{1}}$,
called a \textbf{string operator}.

A finite set $\Omega$ of index vectors is called \textbf{fit} if
for each $i\in\{1,2,\dots,m\}$, there exists a vector $t=(t_{1},t_{2},\dots,t_{q})\in\Omega$
such that $t_{s}=i$ for some $s\in\{1,2,\dots,q\}$. Denote by $\mathcal{M}$
the collection of all pairs $(\Omega,w)$, where $\Omega$ is a finite
fit set of index vectors and $w:\Omega\rightarrow(0,\infty)$ is such
that $\sum_{t\in\Omega}w(t)=1.$

For any $(\Omega,w)\in\mathcal{M}$ define the convex combination
of the end-points of all strings defined by members of $\Omega$
\begin{equation}
P_{\Omega,w}(x):=\sum_{t\in\Omega}w(t)P[t](x),\;x\in X.
\end{equation}

Let $\Delta\in(0,1/m)$ and an integer $\bar{q}\geq m$ be arbitrary
fixed and denote by $\mathcal{M}_{\ast}\equiv\mathcal{M}_{\ast}(\Delta,\bar{q})$
the set of all $(\Omega,w)\in\mathcal{M}$ such that the lengths of
the strings are bounded and the weights are all bounded away from
zero, i.e.,
\begin{equation}
\mathcal{M}_{\ast}:=\{(\Omega,w)\in\mathcal{M\mid}\text{ }\ell(t)\leq\bar{q}\text{ and }w(t)\geq\Delta,\text{ }\forall\text{ }t\in\Omega\}.
\end{equation}

\begin{algorithm}
\label{alg:DSAP}$\left.{}\right.$\textbf{Algorithm 2. The DSAP method
with variable strings and variable weights}

\textbf{Initialization}: Select an arbitrary $x^{0}\in X$,

\textbf{Iterative step}: Given a current iteration vector $x^{k}$
pick a pair $(\Omega_{k},w_{k})\in\mathcal{M}_{\ast}$ and calculate
the next iteration vector $x^{k+1}$ by
\begin{equation}
x^{k+1}=P_{\Omega_{k},w_{k}}(x^{k})\text{.}
\end{equation}
\end{algorithm}

The first prototypical string-averaging algorithmic scheme appeared
in \cite{ceh01} and subsequent work on its realization with various
algorithmic operators includes \cite{CS08,CS09,ct03,cz-2-2014,crombez,gordon,pen09,pscr10,rhee03},
see also \cite{bargetz-dsap-2018} and Zaslavski's book \cite{zaslavski-book-SA}.
If in the DSAP method one uses only a single index vector $t=(1,2,\dots,m)$
that includes all constraints indices then the fully-sequential Kaczmarz
cyclic projection method is obtained. For linear hyperplanes as constraints
sets the latter is equivalent with the, independently discovered,
ART (for Algebraic Reconstruction Technique) in image reconstruction
from projections, see \cite{GTH}. If, at the other extreme, one uses
exactly $m$ index vectors $t=(i),$ for $i=1,2,\dots,m,$ each consisting
of exactly one constraint index, then the fully-simultaneous projection
method of Cimmino is recovered. In-between these ``extremes'' the
DSAP method allows for a large arsenal of specific feasibility-seeking
projection algorithms. See \cite{bb96,annotated,cccdh10} for more
information on projection methods.

The \textbf{superiorized version of the DSAP algorithm} is obtained
by using Algorithm 2 as the algorithmic operator $\mathcal{A}$ in
Algorithm 1. The following result about its behavior was proved. Consider
the set $C_{min}:=\{x\in C\mid\;\phi(x)\leq\phi(y){\text{ for all }}y\in C\},$
and assume that $C_{min}\not=\emptyset.$

\begin{theorem} \label{thm:fejer}\cite[Theorem 4.1]{cz14-feje}
Let $\phi:X\rightarrow R$ be a convex continuous function, and let
$C_{\ast}\subseteq C_{min}$ be a nonempty subset. Let $r_{0}\in(0,1]$
and $\bar{L}\geq1$ be such that, ${\text{for all }}x\in C_{\ast}{\text{ and all }}y$
such that$\;||x-y||\leq r_{0},$
\begin{equation}
|\phi(x)-\phi(y)|\leq\bar{L}||x-y||{,}
\end{equation}
and suppose that $\{(\Omega_{k},w_{k})\}_{k=0}^{\infty}\subset\mathcal{M}_{\ast}.$
Then any sequence $\{y^{k}\}_{k=0}^{\infty},$ generated by the superiorized
version of the DSAP algorithm, converges in the norm of $X$ to a
$y^{\ast}\in C$ and exactly one of the following two alternatives
holds:

(a) $y^{\ast}\in C_{min}$;

(b) $y^{\ast}\notin C_{min}$ and there exist a natural number $k_{0}$
and a $c_{0}\in(0,1)$ such that for each $x\in C_{\ast}$ and for
each integer $k\geq k_{0}$,
\begin{equation}
\Vert y^{k+1}-x\Vert^{2}\leq\Vert y^{k}-x\Vert^{2}-c_{0}\sum_{n=1}^{N_{k}-1}\beta_{k,n}.
\end{equation}

\end{theorem}

This shows that $\{y^{k}\}_{k=0}^{\infty}$ is strictly Fej\'er-monotone
with respect to\textbf{ }$C_{\ast},$ i.e., that $\Vert y^{k+1}-x\Vert^{2}<\Vert y^{k}-x\Vert^{2},$
for all $k\geq k_{0},$ because $c_{0}\sum_{n=1}^{N_{k}-1}\beta_{k,n}>0.$
The strict Fej\'er-monotonicity however does not guarantee convergence
to a constrained minimum point but only says that the so-created feasibility-seeking
sequence $\{y^{k}\}_{k=0}^{\infty}$ has the additional property of
getting strictly closer, without necessarily converging, to the points
of a subset of the solution set of of the constrained minimization
problem.

Published experimental results repeatedly confirm that global reduction
of the value of the objective function $\phi$ is indeed achieved,
without loosing the convergence toward feasibility, see \cite{bk13,bdhk07,cdh10,cz12,dhc09,rand-conmath,gh13,hd08,hgdc12,wenma13,ndh12,pscr10}.
In some of these cases the SM returns a lower value of the objective
function $\phi$ than an exact minimization method with which it is
compared, e.g., \cite{cdhst14}.

\subsection{Strong superiorization\label{sect:strongSM}}

In this section we present a restricted version of the SM of \cite{hgdc12}
as adapted to the situation in \cite{cdhst14}. We consider again
Case \ref{case:cfp} in Section \ref{sec:sup-framework} wherein $T$
is a convex feasibility problem. Let $\varTheta:=\{C_{i}\}{}_{i=1}^{m}$
be a family of nonempty closed convex subsets of a Hilbert space $X$
and set $C=\cap_{i=1}^{m}C_{i}$. We do not assume that $C\neq\emptyset,$
but only that there is some nonempty subset $\Lambda\in R^{n}$ such
that $C_{i}\subseteq\Lambda$ for all $i.$ Instead of the nonemptiness
assumption we associate with the family $\{C_{i}\}_{i=1}^{m}$ a \textbf{proximity
function} ${Prox}_{\varTheta}:\Lambda\rightarrow\mathbb{R}_{+}$ that
is an indicator of how incompatible an $x\in\Lambda$ is with the
constraints. For any given $\varepsilon>0$, a point $x\in\Lambda$
for which ${Prox}_{\varTheta}(x)\leq\varepsilon$ is called an \textbf{$\varepsilon$-compatible
solution for $\varTheta$}. We further assume that we have a feasibility-seeking
algorithmic operator $\mathcal{A}:R^{n}\rightarrow\Lambda$, with
which we define the basic algorithm as the iterative process
\begin{equation}
x^{k+1}=\mathcal{A}(x^{k}),\text{ for all }k\geq0,\text{ for an arbitrary }x^{0}\in\Lambda.
\end{equation}
The following definition helps to evaluate the output of the basic
algorithm upon termination by a stopping rule.

\begin{definition} \textbf{\label{def:epsilon-output}The }$\varepsilon$\textbf{-output
of a sequence. }Given $\varTheta$ such that $C_{i}\subseteq\Lambda\subseteq R^{n}$,
for all $i,$ a proximity function ${Prox}_{\varTheta}:\Lambda\rightarrow R_{+}$,
a sequence $\left\{ x^{k}\right\} _{k=0}^{\infty}\subset\Lambda$
and an $\varepsilon>0,$ then an element $x^{K}$ of the sequence
which has the properties: (i) ${Prox}_{\varTheta}\left(x^{K}\right)\leq\varepsilon,$
and (ii) ${Prox}_{\varTheta}\left(x^{k}\right)>\varepsilon$ for all
$0\leq k<K,$ is called an \textbf{$\varepsilon$-output of the sequence
$\left\{ x^{k}\right\} _{k=0}^{\infty}$ with respect to the pair
$(\varTheta,$ ${Prox}_{\varTheta})$.} \end{definition}

We denote the \textbf{$\varepsilon$-output} by $O\left(\varTheta,\varepsilon,\left\{ x^{k}\right\} _{k=0}^{\infty}\right)=x^{K}.$
Clearly, an $\varepsilon$-output $O\left(\varTheta,\varepsilon,\left\{ x^{k}\right\} _{k=0}^{\infty}\right)$
of a sequence $\left\{ x^{k}\right\} _{k=0}^{\infty}$ might or might
not exist, but if it does, then it is unique. If $\left\{ x^{k}\right\} _{k=0}^{\infty}$
is produced by an algorithm intended for the feasible set $C,$ such
as the Basic Algorithm, without a termination criterion, then $O\left(\varTheta,\varepsilon,\left\{ x^{k}\right\} _{k=0}^{\infty}\right)$
is the output produced by that algorithm when it includes the termination
rule to stop when an $\varepsilon$-compatible solution for $\varTheta$
is reached.

\begin{definition} \textbf{Strong perturbation resilience. }Given
$\varTheta$ such that $C_{i}\subseteq\Lambda\subseteq R^{n}$, for
all $i,$ a proximity function ${Prox}_{\varTheta}:\Lambda\rightarrow R_{+}$,
an algorithmic operator $\mathcal{A}$ and an $x^{0}\in\Lambda$.
We use $\left\{ x^{k}\right\} _{k=0}^{\infty}$ to denote the sequence
generated by the Basic Algorithm when it is initialized by $x^{0}$.
The basic algorithm is said to be\texttt{ }\textbf{strongly perturbation
resilient} iff the following hold: (i) there exist an $\varepsilon>0$
such that the $\varepsilon$-output $O\left(\varTheta,\varepsilon,\left\{ x^{k}\right\} _{k=0}^{\infty}\right)$
exists for every $x^{0}\in\Lambda$; (ii) for every $\varepsilon>0,$
for which the $\varepsilon$-output $O\left(\varTheta,\varepsilon,\left\{ x^{k}\right\} _{k=0}^{\infty}\right)$
exists for every $x^{0}\in\Lambda$ it holds that the $\varepsilon^{\prime}$-output
$O\left(\varTheta,\varepsilon^{\prime},\left\{ y^{k}\right\} _{k=0}^{\infty}\right)$
exists for every $\varepsilon^{\prime}>\varepsilon$ and for every
sequence $\left\{ y^{k}\right\} _{k=0}^{\infty}$ generated by
\begin{equation}
y^{k+1}=\mathcal{A}\left(y^{k}+\beta_{k}v^{k}\right),\text{ for all }k\geq0,\label{eq:perturb}
\end{equation}
where the vector sequence $\left\{ v^{k}\right\} _{k=0}^{\infty}$
is bounded and the scalars $\left\{ \beta_{k}\right\} _{k=0}^{\infty}$
are such that $\beta_{k}\geq0$, for all $k\geq0,$ and $\sum_{k=0}^{\infty}\beta_{k}<\infty$.
\\
\end{definition}

A theorem which gives sufficient conditions for strong perturbation
resilience of the basic algorithm has been proved in \cite[Theorem 1]{hgdc12}. 

Along with $\varTheta$ such that $C_{i}\subseteq\Lambda\subseteq R^{n}$,
for all $i$ and a proximity function ${Prox}_{\varTheta}:\Lambda\rightarrow R_{+}$,
we look at the objective function $\phi:R^{n}\rightarrow R$, with
the convention that a point in $R^{n}$ whose value of $\phi$ is
smaller is considered \textbf{superior} to a point in $R^{n}$ for
which the value of $\phi$ is larger. The essential idea of the SM
is to make use of the perturbations of (\ref{eq:perturb}) to transform
a strongly perturbation resilient Basic Algorithm that seeks a constraints-compatible
solution for $\varTheta$ into its Superiorized Version whose outputs
are equally good from the point of view of constraints-compatibility,
but are superior (not necessarily optimal) according to the objective
function $\phi$.

\begin{definition} \label{def:nonascend}Given a function $\phi:R^{n}\rightarrow R$
and a point $y\in R^{n}$, we say that a vector $d\in R^{n}$ is \textbf{nonascending
for $\phi$ at}\texttt{ }$y$ iff $\left\Vert d\right\Vert \leq1$
and there is a $\delta>0$ such that for all $\lambda\in\left[0,\delta\right]$
we have $\phi\left(y+\lambda d\right)\leq\phi\left(y\right).$\\
\end{definition}

Obviously, the zero vector is always such a vector, but for superiorization
to work we need a sharp inequality to occur in (\ref{def:nonascend})
frequently enough. The Superiorized Version of the Basic Algorithm
assumes that we have available a summable sequence $\left\{ \eta_{\ell}\right\} _{\ell=0}^{\infty}$
of positive real numbers (for example, $\eta_{\ell}=a^{\ell}$, where
$0<a<1$) and it generates, simultaneously with the sequence $\left\{ y^{k}\right\} _{k=0}^{\infty}$
in $\Lambda$, sequences $\left\{ v^{k}\right\} _{k=0}^{\infty}$
and $\left\{ \beta_{k}\right\} _{k=0}^{\infty}$. The latter is generated
as a subsequence of $\left\{ \eta_{\ell}\right\} _{\ell=0}^{\infty}$,
resulting in a nonnegative summable sequence $\left\{ \beta_{k}\right\} _{k=0}^{\infty}$.
The algorithm further depends on a specified initial point $y^{0}\in\Lambda$
and on a positive integer $N$. It makes use of a logical variable
called \textit{loop}\emph{. }The general form of the superiorized
version of the basic algorithm is presented next by its pseudo-code.

\begin{algorithm}
\label{alg_super}\textbf{Algorithm 3. General form of the superiorized
version of the basic algorithm} 
\end{algorithm}

\begin{enumerate}
\item \textbf{set} $k=0$
\item \textbf{set} $y^{k}=y^{0}$
\item \textbf{set} $\ell=-1$
\item \textbf{repeat}
\item $\qquad$\textbf{set} $n=0$
\item $\qquad$\textbf{set} $y^{k,n}=y^{k}$
\item $\qquad$\textbf{while }$n$\textbf{$<$}$N$
\item $\qquad$\textbf{$\qquad$set }$v^{k,n}$\textbf{ }to be a nonascending
vector for $\phi$ at $y^{k,n}$
\item $\qquad$\textbf{$\qquad$set} \emph{loop=true}
\item $\qquad$\textbf{$\qquad$while}\emph{ loop}
\item $\qquad\qquad\qquad$\textbf{set $\ell=\ell+1$}
\item $\qquad\qquad\qquad$\textbf{set} $\beta_{k,n}=\eta_{\ell}$
\item $\qquad\qquad\qquad$\textbf{set} $z=y^{k,n}+\beta_{k,n}v^{k,n}$
\item $\qquad\qquad\qquad$\textbf{if }$\phi\left(z\right)$\textbf{$\leq$}$\phi\left(y^{k}\right)$\textbf{
then }
\item $\qquad\qquad\qquad\qquad$\textbf{set }$n$\textbf{$=$}$n+1$
\item $\qquad\qquad\qquad\qquad$\textbf{set }$y^{k,n}$\textbf{$=$}$z$
\item $\qquad\qquad\qquad\qquad$\textbf{set }\emph{loop = false}
\item $\qquad$\textbf{set }$y^{k+1}$\textbf{$=$}$\mathcal{A}\left(y^{k,N}\right)$
\item $\qquad$\textbf{set }$k=k+1$ 
\end{enumerate}
\begin{theorem} \label{theorem4.5-1} Any sequence $\left\{ y^{k}\right\} _{k=0}^{\infty}$,
generated by the superiorized version of the basic algorithm, Algorithm
3, satisfies (\ref{eq:perturb}). Further, if, for a given $\varepsilon>0,$
the $\varepsilon$-output $O\left(\varTheta,\varepsilon,\left\{ x^{k}\right\} _{k=0}^{\infty}\right)$
of the Basic Algorithm exists for every $x^{0}\in\Lambda$, then every
sequence $\left\{ y^{k}\right\} _{k=0}^{\infty}$, generated by the
Algorithm 3, has an $\varepsilon^{\prime}$-output $O\left(\varTheta,\varepsilon^{\prime},\left\{ y^{k}\right\} _{k=0}^{\infty}\right)$
for every $\varepsilon^{\prime}>\varepsilon$. \end{theorem}

The proof of this theorem follows from the analysis of the behavior
of the superiorized version of the basic algorithm in \cite[pp. 5537--5538]{hgdc12}.
In other words, Algorithm 3 produces outputs that are essentially
as constraints-compatible as those produced by the original basic
algorithm. However, due to the repeated steering of the process by
lines 7 to 17 toward reducing the value of the objective function
$\phi$, we can expect that its output will be superior (from the
point of view of $\phi$) to the output of the (unperturbed) Basic
Algorithm.

Algorithm 1 and Algorithm 3 are not identical but are based on the
same leading principle of the superiorization methodology. Comments
on the differences between them can be found in \cite[Remark 4.1]{cz14-feje}.
Nevertheless, the Theorems \ref{thm:fejer} and \ref{theorem4.5-1}
related to these superiorized versions of the basic algorithm, respectively,
leave the question of rigorously analyzing the behavior of the SM,
under various conditions, open.

\subsection{Controlling the effect of the perturbations\label{subsec:Controlling}}

The scalars $\beta_{k}$ in the SM algorithmic scheme, see Definition
\ref{def:resilient} and Eq. (\ref{eq:perturb}), are generated such
that $\beta_{k}\geq0$, for all $k\geq0,$ and $\sum_{k=0}^{\infty}\beta_{k}<\infty$.
This implies that they form a tending to zero sequence. As step-sizes
of the perturbation the effects of the objective function decrease
is bound to diminish as iterations proceed.

In some applications, various methods have been studied for controlling
the step-sizes, see, e.g., \cite{Langthaler,Prommegger}, see also
the software package SNARK14 \cite{SNARK} which is an updated version
of \cite{SNARK-09}. Recently, a new strategy which allows restarting
the sequence of step-sizes to a previous value while maintaining the
summability of the series of step-sizes was suggested \cite{restarts},
resulting in improvement of the algorithm\textquoteright s performance.

\section{Derivative-free superiorization}

\subsection{Derivative-free superiorization and derivative-free optimization}

Here we describe the general applicability of \textbf{derivative-free
superiorization} (DFS) as an alternative to previously proposed superiorization
approaches. These earlier approaches were based on generation of nonascending
vectors, for objective function reduction steps, that mostly required
the ability to calculate gradients or subgradients of the objective
function. Observing the body of knowledge of \textbf{derivative-free
optimization} (DFO), see, e.g., \cite{Conn-book-2009}, we explore
a DFS algorithm.

In DFS, the perturbation phase of the superiorized version of a basic
algorithm the objective function reduction steps that depend on gradient
or subgradient calculations are replaced by steps that use a direction
search technique which does not require any form of differentiability.
Continuing the work of \cite{CHS18}, we searched in \cite{CGHH-DFS}
the neighborhood of a current point $x$ for a point at which the
objective function exhibits nonascent.

While this might seem a simple technical matter, the ramifications
for practical applications of the SM are important. For example, in
intensity-modulated radiation therapy treatment planning, with photons,
protons or other particles, the \textbf{normal tissue complication
probability} (NTCP) is a predictor of radiobiological effects for
organs at risk. The inclusion of it, or of other biological functions,
as an objective function in the mathematical problem modeling and
the planning algorithm, is hampered because they are, in general,
empirical functions whose derivatives cannot be calculated, see, e.g.,
\cite{gay-niemierko-ntcp-2007}. In the recent paper \cite{nystorm-2020}
the authors list issues of immediate clinical and practical relevance
to the Proton Therapy community, highlighting the needs for the near
future but also in a longer perspective. They say that ``...practical
tools to handle the variable biological efficiency in Proton Therapy
are urgently demanded...''.

The output of a superiorized version of a constraints-compatibility-seeking
algorithm will have smaller (but not minimal) objective function $\phi$
value than the output of the same constraints-compatibility-seeking
algorithm without perturbations, everything else being equal. Even
though superiorization is not an exact minimization method, we think
of it as an applicable (and possibly, more efficacious) alternative
to derivative-free constrained minimization methods applied to the
same data for two main reasons: its ability to handle constraints
and its ability to cope with very large-size problems. This is in
contrast with the current state of the art, which is as follows.

The review paper of Rios and Sahinidis \cite{Rios} ``... addresses
the solution of \emph{bound-constrained} optimization problems using
algorithms that require only the availability of objective function
values but no derivative information,'' with bound constraints imposed
on the vector $x$. The book by Conn, Scheinberg and Vicente \cite{Conn-book-2009}
deals only with derivative-free unconstrained minimization, except
for its last chapter (of 10 pages out of the 275) entitled ``Review
of constrained and other extensions to derivative-free optimization.''
Li \emph{et al}. \cite{LCLLLL} do not even mention constraints. In
\cite{diniz2011} the numerical work deals with: ``The dimension
of the problems {[}i.e., the size of the vector $x${]} varies between
2 and 16, while the number of constraints are between 1 and 38, exceeding
10 in only 5 cases.'' In \cite{dfo-4-oil} the numerical tests are
limited to: ``The first case has 80 optimization variables {[}i.e.,
the size of the vector $x${]} and only bound constraints, while the
second example is a generally constrained production optimization
involving 20 optimization variables and 5 general constraints.''
Similar orders of magnitude for problem sizes appear in the numerical
results presented in \cite{Audet-Dennis-2009} and also in the book
of Audet and Hare \cite{Audet-book-2017}.

This indicates that (i) much of the literature on derivative-free
minimization is concerned with unconstrained minimization or with
bound-constraints on the variables, and (ii) many, if not all, proposed
methods were designed (or, at least, demonstrated) only for small-scale
problems. In contrast, the DFS method proposed here can handle any
type of constraints for which a separate efficient constraints-compatibility-seeking
algorithm is available and is capable of solving very large problems.
In the matter of problem sizes, we discover here, admittedly with
a very preliminary demonstration, that DFS can compete well with DFO
on large problems. Since the constraints-compatibility-seeking algorithm
forms part of the proposed DFS method, the method can use exterior
initialization (that is initializing the iterations at any point in
space). Furthermore, very large-scale problems can be accommodated.

The \textbf{progressive barrier} (PB) approach, described in Chapter
12 of the book \cite{Audet-book-2017}, originally published in \cite{Audet-Dennis-2009},
is an alternative to the \textbf{exterior penalty} (EP) approach that
is mentioned in \cite{CGHH-DFS}. However, the PB differs from our
DFS method, in spite of some similarities with it, as explained in
\cite{CGHH-DFS}.

\subsection{\label{sec:proximity-target-curve}The proximity-target curve}

A tool for deciding which of two iterative methods is ``better''
for solving a particular problem was presented in \cite{CGHH-DFS}
and applied to the DFS algorithm developed there. We care to reproduce
it here because of its potential usefulness in other situations where
two algorithms are compared. Since an iterative method produces a
sequence of points, our tool is based on such sequences. Furthermore,
since in the SM we are interested in the values of two functions (proximity
function and target function) at each iteration point, the efficacy
of the behavior of the iterative method can be represented by a curve
in two-dimensional space, defined as \textbf{the proximity-target
curve} below. 

It indicates the target value for any achieved proximity value. This
leads to the intuitive concept of an algorithm being ``better''
than another one, if its proximity target curve is below that of the
other one (that is, the target value for it is always smaller than
the target value of the other one for the same proximity value). Such
may not always be the case, the two proximity curves may cross each
other, providing us with intervals of proximity values within which
one or the other method is better.

While this way of thinking is not common in numerical analysis and
in optimization, it is quite generally used in many sciences in situations
where it is desirable to obtain an object for which the values of
two evaluating functions are small simultaneously. A prime example
is in estimation theory where we desire an estimation method with
both small bias and small variance. More specifically, is the concept
of a \textbf{receiver operating characteristics (ROC) curve} that
illustrates the diagnostic ability of a binary classifier system as
its discrimination threshold is varied. It is created by plotting
the true positive rate against the false positive rate at various
threshold settings. One classifier system is considered ``better''
than the other one if its ROC curve is above that of the other one;
but, just as for our proximity-target curves, the ROC curves for two
classifier systems may cross each other. There are many publications
on the role of ROC curves in the evaluation of medical imaging techniques;
see, for example, \cite{METZ89a,SWET79a}. Their use for image reconstruction
algorithm evaluation is discussed, for example, in \cite{COOL92a}.

\begin{figure}
\begin{centering}
\includegraphics[scale=0.1]{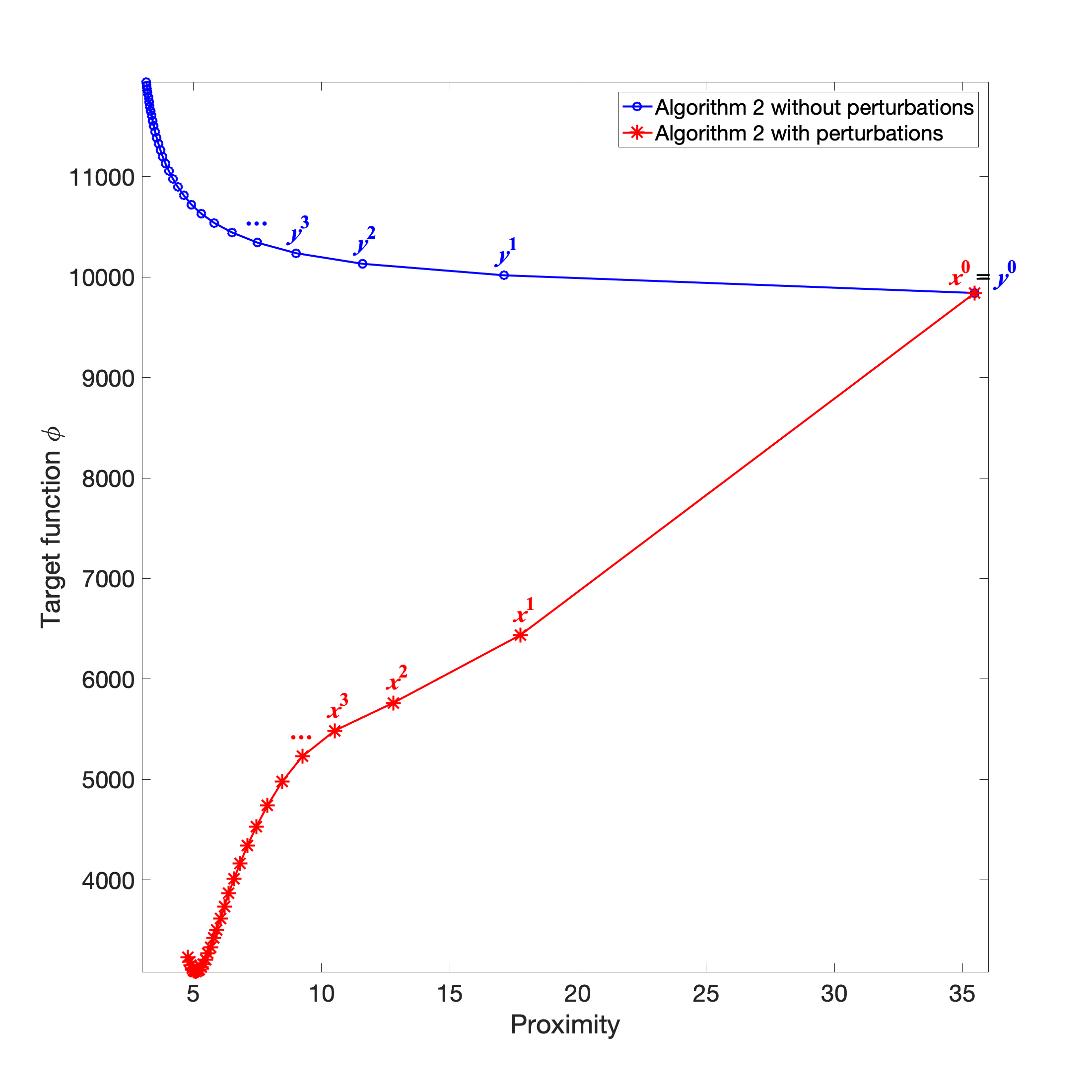}
\par\end{centering}
\caption{\label{fig:Target-function-values}Proximity-target curves $P$ and
$Q$ of the first 30 iterates of a feasibility-seeking algorithm (Algorithm
2 in \cite{CGHH-DFS}) with perturbations (\EightStarTaper ) and without
perturbations ($\circ)$. Reproduced from \cite{CGHH-DFS}.}
\end{figure}

For incarnations of the definitions the reader may wish to look ahead
to Figure \ref{fig:Target-function-values}. That figure illustrates
the discussed notions for two particular finite sequences $U:=\left(x^{k}\right)_{k=K_{lo}}^{K_{hi}}$
and $V:=\left(y^{k}\right)_{k=L_{lo}}^{L_{hi}}$. The details of how
those sequences were specifically produced are given below in \cite{CGHH-DFS}.
We use the notations of Subsection \ref{sect:strongSM} and adapt
the, more general, definitions of \cite{CGHH-DFS} to this case.

\begin{definition}\textbf{Monotone proximity of a finite sequence}.
Consider $\varTheta$ such that $C_{i}\subseteq\Lambda\subseteq R^{n}$
for all $i,$ and a proximity function ${Prox}_{\varTheta}:\Lambda\rightarrow R_{+}$.
For positive integers $K_{lo}$ and $K_{hi}>K_{lo}$, the finite sequence
$U:=\left(x^{k}\right)_{k=K_{lo}}^{K_{hi}}$ of points in $\Lambda$
is said to be \textbf{of monotone proximity} if for $K_{lo}<k\leq K_{hi}$,
${Prox}_{\varTheta}\left(x^{k-1}\right)>{Prox}_{\varTheta}\left(x^{k}\right)$.\\
\end{definition}
\begin{definition}\textbf{The proximity-target curve of a finite
sequence}. Consider $\varTheta$ such that $C_{i}\subseteq\Lambda\subseteq R^{n}$
for all $i,$ a proximity function ${Prox}_{\varTheta}:\Lambda\rightarrow R_{+}$,
a target function $\phi:\Lambda\rightarrow R$, and positive integers
$K_{lo}$ and $K_{hi}>K_{lo}$. Let $U:=\left(x^{k}\right)_{k=K_{lo}}^{K_{hi}}$
be a sequence of monotone proximity. Then the \textbf{proximity-target
curve} $P\subseteq R^{2}$ associated with $U$ is uniquely defined
by:
\begin{enumerate}
\item For $K_{lo}\leq k\leq K_{hi}$, $\left({Prox}_{\varTheta}\left(x^{k}\right)\!,\phi\left(x^{k}\right)\right)\in P$.
\item The intersection $\{(y,x)\in R^{2}\mid{Prox}_{\varTheta}\left(x^{k}\right)\leq y\leq{Prox}_{\varTheta}\left(x^{k-1}\right)\}\cap P$
is the line segment from $\left({Prox}_{\varTheta}\left(x^{k-1}\right),\phi\left(x^{k-1}\right)\right)$
\end{enumerate}
to $\left({Prox}_{\varTheta}\left(x^{k}\right),\phi\left(x^{k}\right)\right)$.\\
\end{definition}
to $\left({Prox}_{\varTheta}\left(x^{k}\right),\phi\left(x^{k}\right)\right)$.\\
\begin{definition}\label{def:comparison}\textbf{Comparison of proximity-target
curves of sequences}

Consider $\varTheta$ such that $C_{i}\subseteq\Lambda\subseteq R^{n}$
for all $i,$ a proximity function ${Prox}_{\varTheta}:\Lambda\rightarrow R_{+}$,
a target function $\phi:\Lambda\rightarrow R$, and positive integers
$K_{lo}$, $K_{hi}>K_{lo}$, $L_{lo}$, $L_{hi}>L_{lo}$, let $R:=\left(x^{k}\right)_{k=K_{lo}}^{K_{hi}}$
and $S:=\left(y^{k}\right)_{k=L_{lo}}^{L_{hi}}$ be sequences of points
in $\Omega$ of monotone proximity for which $P$ and $Q$ are their
respective associated proximity-target curves. Define
\begin{equation}
\begin{array}{c}
t:=\max\left({Prox}_{\varTheta}\left(x^{K_{hi}}\right),{Prox}_{\varTheta}\left(y^{L_{hi}}\right)\right),\\
u:=\min\left({Prox}_{\varTheta}\left(x^{K_{lo}}\right),{Prox}_{\varTheta}\left(y^{L_{lo}}\right)\right).
\end{array}\label{eq:limits-1}
\end{equation}
Then \textbf{$R$ is better targeted than $S$} if:
\begin{enumerate}
\item $t\leq u$ and
\item for any real number $h$, if $t\leq h\leq u$, $\left(h,v\right)\in P$
and $\left(h,w\right)\in Q$, then $v\leq w$.
\end{enumerate}
\end{definition}

\medskip{}
This definition is intuitively desirable. Suppose that we have an
iterative algorithm that produces a sequence, $y^{0},y^{1},y^{2},\cdots$,
of which $S:=\left(y^{k}\right)_{k=L_{lo}}^{L_{hi}}$ is a subsequence.
An alternative algorithm that produces a sequence of points of which
$R:=\left(x^{k}\right)_{k=K_{lo}}^{K_{hi}}$ is a subsequence that
is better targeted than $S$ has a desirable property: Within the
range $\left[t,u\right]$ of proximity values, the point that is produced
by the alternative algorithm with that proximity value, is likely
to have a lower (and definitely not higher) value of the target function
than the point with that proximity value that is produced by the original
algorithm. This property is stronger than what we stated before, namely
that superiorization produces an output that is equally good from
the point of view of proximity, but is superior with respect to the
target function. Here the single output determined by a fixed $\varepsilon$
is replaced by a set of potential outputs for any $\varepsilon\in\left[t,u\right]$.

\section{The superiorization methodology guarantee problem\label{sect:frame-1}}

It is clear from the above text that the SM is interlacing the iterative
steps of two separate and independent iterative processes. Therefore,
we reformulate here \textbf{the guarantee problem of the SM} in the
following general terms. One iterative process, the basic algorithm,
is of the form
\begin{equation}
x^{0}\in R^{n},\;x^{k+1}=\mathcal{A}(x^{k}),\;k\geq0,\label{eq:basic-alg-1}
\end{equation}
where $\mathcal{A}:R^{n}\rightarrow R^{n}$ is a given algorithmic
operator and $x^{0}\in R^{n}$ is an arbitrary initialization point.
It is assumed that this is a convergent process
\begin{equation}
\lim_{k\rightarrow\infty}x^{k}=x^{\ast}\in S\subseteq R^{n},
\end{equation}
where $S\neq\varnothing$ and $S$ is not a singleton. $S$ could
be a solution set of some problem, e.g., the convex feasibility problem.
It is further assumed that the basic algorithm is bounded perturbations
resilient as in Definition \ref{def:resilient}.

The second iterative process, henceforth called ``the auxiliary algorithm'',
is 

\begin{equation}
z^{0}\in R^{n},\;z^{k+1}=\mathcal{B}(z^{k}),\;k\geq0,\label{eq:second-alg-1}
\end{equation}
where $\mathcal{B}:R^{n}\rightarrow R^{n}$ is another given algorithmic
operator and $z^{0}\in R^{n}$ is an arbitrary initialization point.
It is assumed that this is also a convergent process
\begin{equation}
\lim_{k\rightarrow\infty}z^{k}=z^{\ast}\in T\subseteq R^{n},
\end{equation}
where $T\neq\varnothing$ and $T$ is not a singleton. $T$ could
be a solution set of some other problem.

In these circumstances a general description of the SM consists of
perturbing the iterates of the basic algorithm (\ref{eq:basic-alg-1})
by perturbations induced by the auxiliary algorithm (\ref{eq:second-alg-1}),
i.e.,
\begin{equation}
v^{k}:=\frac{\mathcal{B}(y^{k})-y^{k}}{\left\Vert \mathcal{B}(y^{k})-y^{k}\right\Vert },\;\textup{if}\;\left\Vert \mathcal{B}(y^{k})-y^{k}\right\Vert \neq0,\;\textup{and}\;v^{k}:=0,\;\textup{otherwise}.\label{eq:interlaced}
\end{equation}
This yields the superiorized version of the basic algorithm

\begin{equation}
y^{k+1}=\mathcal{A}\left(y^{k}+\beta_{k}v^{k}\right),\text{ for all }k\geq0,\label{eq:superiorized-general}
\end{equation}
with $v^{k}$ as in (\ref{eq:interlaced}) and$\left\{ \beta_{k}\right\} _{k=0}^{\infty}$
such that $\beta_{k}\geq0$, for all $k\geq0,$ and $\sum_{k=0}^{\infty}\beta_{k}<\infty$.

\textbf{Example}. If the operator $\mathcal{A}$ is a feasibility-seeking
algorithm for a given convex feasibility problem and $\mathcal{B}$
is an unconstrained gradient descent algorithm for a given objective
function $f:R^{n}\rightarrow R$ then (\ref{eq:superiorized-general})
describes the earlier presented SM algorithms.

In \cite[Section 3]{censor-levy-2019} we gave a precise definition
of \textbf{the guarantee problem of the SM}. We wrote there: ``The
SM interlaces into a feasibility-seeking basic algorithm target function
reduction steps. These steps cause the target function to reach lower
values locally, prior to performing the next feasibility-seeking iterations.
A mathematical guarantee has not been found to date that the overall
process of the superiorized version of the basic algorithm will not
only retain its feasibility-seeking nature but also accumulate and
preserve globally the target function reductions.'' 

\begin{definition}\textbf{The guarantee problem of the SM}. Under
which conditions one can guarantee that a superiorized version of
a bounded perturbation resilient feasibility-seeking algorithm converges
to a feasible point that has target function value smaller or equal
to that of a point to which this algorithm would have converged if
no perturbations were applied -- everything else being equal.

\end{definition}

Numerous works that are cited in \cite{SM-bib-page} show that this
global function reduction of the SM occurs in practice in many real-world
applications. In addition to a partial answer in \cite{censor-levy-2019}
with the aid of \textbf{the concentration of measure principle} there
is also the partial result of \cite[Theorem 4.1]{cz14-feje} about
\textbf{strict Fej\'er monotonicity} of sequences generated by an SM
algorithm.

Proving mathematically a guarantee of global function reduction of
the SM will probably require some additional assumptions on the feasible
set, on the objective function, on the parameters involved, or even
on the set of permissible initialization points.

\section{Some applications of superiorization\label{sec:applications}}

The, compiled and continuously updated, bibliography \cite{SM-bib-page}
contains references to numerous about superiorization (161 items as
of December 4, 2022). We single out a few that applied the SM successfully
to significant real-world problems. For many more theoretical, experimental
studies or papers dealing with a variety of real-world applications,
we recommend \cite{SM-bib-page}.

In \cite{guenter-compare-sup-reg-2022} Guenter et al. consider the
fully-discretized modeling of \textbf{image reconstruction from projections}
problem that leads to a system of linear equations which is huge and
very sparse. Solving such systems, sometimes under limitations on
the computing resources, has been, is, and will remain a challenge.
The authors aim not only at solving the linear system resulting from
the modeling alone but consider the constrained optimization problem
of minimizing an objective function subject to the modeling constraints.
To do so, they recognize two fundamental approaches: (i) superiorization,
and (ii) regularization. Within these two methodological approaches
they evaluate 21 algorithms over a collection of 18 different phantoms
(i.e., test problems), presenting their experimental results in very
informative ways.

In \cite{fink-2021} Fink et al. study the \textbf{nonconvex multi-group
multicast beamforming problem} with quality-of-service constraints
and per-antenna power constraints. They formulate a convex relaxation
of the problem as a semidefinite program in a real Hilbert space,
which allows them to approximate a point in the feasible set by iteratively
applying a bounded perturbation resilient fixed-point mapping. Inspired
by the superiorization methodology, they use this mapping as a basic
algorithm, and add in each iteration a small perturbation with the
intent to reduce the objective value and the distance to nonconvex
rank-constraint sets.

Pakkaranang et al. \cite{Pakkaranang-2020} construct a novel algorithm
for solving \textbf{non-smooth composite optimization problems}. By
using an inertial technique, they propose a modified proximal gradient
algorithm with outer perturbations and obtain strong convergence results
for finding a solution of composite optimization problem. Based on
bounded perturbation resilience, they present their algorithm with
the superiorization method and apply it to image recovery problem.
They provide numerical experiments that show the efficiency of the
algorithm and compare it with previously known algorithms in signal
recovery.

In \cite{barkmann-2022} the SM is applied to the \textbf{intensity-modulated
radiation therapy (IMRT) treatment planning} problem. It is found
there that the superiorization prototype solved the linearly constrained
planning problem with similar performance to that of a general purpose
nonlinear constrained optimizer while showing smooth convergence in
both constraint proximity and objective function reduction. The authors'
work shows that superiorization is a useful alternative to constrained
optimization in radiotherapy inverse treatment planning.

Especially interesting is the recent work of Ma et al. \cite{sahinidis-JOGO-2021}
who propose a novel decomposition framework for \textbf{derivative-free
optimization (DFO) algorithms} which significantly extends the scope
of current DFO solvers to larger-scale problems. They show that their
proposed framework closely relates to the superiorization methodology
that has been used for improving the efficiency of feasibility-seeking
algorithms for constrained optimization problems in a derivative-based
setting. 

\section{Concluding remarks\label{sect: conclusion}}

In many mathematical formulations of significant real-world technological
or physical problems, the objective function is exogenous to the modeling
process which defines the constraints. In such cases, the faith of
the modeler in the usefulness of an objective function for the application
at hand is limited and, as a consequence, it is probably not worthwhile
to invest too much resources in trying to reach an exact constrained
minimum point. This is an argument in favor of using the superiorization
methodology for practical applications. In doing so the amount of
computational efforts invested alternatingly between performing perturbations
and applying the basic algorithm's algorithmic operator can, and needs
to, be carefully controlled in order to allow both activities to properly
influence the outcome. Better theoretical insights into the behavior
of weak and of strong superiorization as well as better ways of implementing
the methodology are needed and await to be developed.

\bigskip{}

\textbf{Acknowledgments.} I am indebted to my partners in the work
on suepriorization, past and present, from whom I have learned so
much. In particular, my thanks go to the collaborators on the papers
on superiorization that are referenced in this paper. Work on this
paper was supported by the ISF-NSFC joint research program Grant No.
2874/19.

\bigskip{}

\end{document}